\documentclass[probth]{svjour1}

\usepackage{amssymb}
\usepackage{amsmath}
\usepackage{amsopn}
\usepackage{mathrsfs}
\usepackage{comment}

\DeclareMathOperator{\supp}{supp} \DeclareMathOperator{\Var}{Var}
\DeclareMathOperator{\diam}{diam}
\newcommand{\norm}[1]{\lVert #1 \rVert}
\newtheorem{assumption}{Assumption}

\begin{document}
  \let\bt\mathbf       
  \let\bm\boldsymbol   
  \let\ca\mathcal      
  \let\scr\mathscr     
  \let\goth\mathfrak   
  \let\bb\mathbb       

  \def\Ex{{\bf E}}
  \def\Pb{{\bf P}}
  \def\Qb{{\bf Q}}
  \def\Norm{{\boldsymbol{\mathcal N}}}
  \def\Rb{{\bf R}}
  \def\KL{{\rm KL}}
  \def\var{{\bf Var}}
  \def\eps{\varepsilon}
  \def\sign{\mathop{\rm sgn}}
  \def\cov{\mathop{\rm Cov}}
  \def\ZZ{\mathbb Z}
  \def\DD{\mathbb D}
  \def\FF{\mathbb F}
  \def\RR{\mathbb R}
  \def\BB{\mathbb B}
  \def\EE{\mathbb E}
  \def\CC{\mathbb C}
  \def\NN{\mathbb N}
  \def\GG{\mathbb G}
  \def\QQ{\mathbb Q}
  \def\HH{\mathbb H}
  \def \1{{\rm 1}\mskip -4,5mu{\rm l} } 
  \def\ds{\displaystyle}
  \def\cal{\mathcal}
  \def\Hbl{\mathcal H(\beta,L)}
  \def\mAb{\mathcal A_b}
  \def\mb{\mathbf b}
  \def\mv{\mathbf v}
  \def\mV{\mathbf V}
  \def\bC{\mathbf C}
  \def\be{\mathbf e}
  \def\bw{\boldsymbol{w}}
  \def\beps{\boldsymbol{\varepsilon}}
  \def\bth{\boldsymbol{\theta}}
  \def\th{\vartheta}
  \def\calT{\mathcal T}
  \def\deg{\mathop{\rm deg}}
  \def\Tr{\mathop{\rm Tr}}

  \setlength{\arraycolsep}{2pt}
  \setlength{\parskip}{.12in}
  \parindent=0pt

  \title{Asymptotic statistical equivalence for ergodic
  diffusions: the multidimensional case}

  \titlerunning{Equivalence for multidimensional diffusions}

  \author{Arnak Dalalyan\inst{1} \and Markus Rei{\ss}\inst{2}}

  \authorrunning{Dalalyan and Rei{\ss}}

\institute{Laboratoire de Probabilit\'es,
  Universit\'e Paris VI,
  Place Jussieu,
  75252 Paris Cedex 05,
  France,
  \email{dalalyan@ccr.jussieu.fr}
\and
  Weierstra{\ss} Institute for Applied Analysis and Stochastics,
  Mohrenstra{\ss}e 39,
  10117 Berlin,
  Germany,
  \email{mreiss@wias-berlin.de}}

  \date{\today}

  \maketitle



\begin{abstract}
  Asymptotic local equivalence in the sense of Le Cam is established for inference on the drift in
  multidimensional ergodic diffusions and an accompanying sequence of Gaussian
  shift experiments. The nonparametric local neighbourhoods can be
  attained for any dimension, provided the regularity of the drift
  is sufficiently large. In addition, a heteroskedastic Gaussian
  regression experiment is given, which is also locally asymptotically
  equivalent and which does not depend on the centre of
  localisation. For one direction of the equivalence an explicit Markov kernel
  is constructed.
\end{abstract}

\section{Introduction}

Asymptotic equivalence is a powerful concept for analysing
statistical inference problems by a transfer to the analogous
problem in a simpler statistical experiment. A breakthrough were the
results by Brown and Low \cite{BrownLow} and Nussbaum
\cite{Nussbaum} who established asymptotic equivalence of the two
classical experiments, one-dimensional Gaussian regression and
density estimation, with an accompanying sequence of Gaussian shift
experiments. In this paper we consider the statistical inference for
the drift in a multidimensional diffusion experiment under
stationarity assumptions and prove the asymptotic equivalence with
corresponding multidimensional Gaussian shift and regression
experiments.

Asymptotic equivalence results for dependent data are not very
numerous, see Dalalyan and Rei{\ss} \cite{DalReiss} for an overview.
Even for simple experiments, as the classical ones described above,
results for asymptotic equivalence in the multidimensional case are
very scarce. We only know of the recent work by Carter \cite{Carter}
who proves asymptotic equivalence for two-dimensional Gaussian
regression, but argues that his method fails for higher dimensions.
One of the main reasons for the difficulties in transferring methods
to higher dimensions is that piecewise constant approximations of
the unknown functional parameter usually do not suffice anymore and
higher order approximations have to be used, which creates
unexpected problems. Brown and Zhang \cite{BrownZhang} remark that
the two classical experiments and their accompanying Gaussian shift
experiments are not asymptotically equivalent in the case of
nonparametric classes of H\"older regularity $\beta\le d/2$, where
$d$ denotes the dimension.

The methodology we applied in \cite{DalReiss} to establish
asymptotic equivalence for scalar diffusions relied heavily on the
concept of local time. For multidimensional diffusions local time
does not exist. This might explain why the statistical theory for
scalar diffusions is very well developed (see Kutoyants
\cite{Kutoyants}), while inference problems for multidimensional
diffusions are more involved and much less studied.  We refer to
Bandi and Moloche \cite{BandiMoloche} for the analysis of kernel
estimators for the drift vector and the diffusion matrix and to
A{\"i}t-Sahalia \cite{AitSahalia} for a recent discussion of
applications for multidimensional diffusion processes in
econometrics.

In Section 2 we review results for multidimensional diffusions and
construct estimators for the invariant density and the drift vector.
Interestingly, the estimator of the invariant density converges for
$d\ge 2$ with a rate which is slower than parametric, but faster
than in classical $d$-dimensional density estimation problems. The
local equivalence result of the multidimensional diffusion
experiment with an accompanying Gaussian shift experiment is
formulated and described in Section 3. The local neighbourhoods can
be attained for drift functions in a nonparametric class of
regularity $\beta>(d-1+\sqrt{2(d-1)^2-1})/2$ for any dimension $d\ge
2$. In Section 4 the corresponding equivalence with a
heteroskedastic regression experiment, which does not depend on the
centre of localisation, is treated. This can be used to establish
global equivalence with a single experiment, which even in the
one-dimensional case cannot be obtained for the Gaussian shift
experiment due to the absence of a variance stabilising transform,
as was first noted by Delattre and Hoffmann \cite{DelHoff}. The
explicit construction of a Markov kernel establishing the important
part of the asymptotic equivalence is presented in Section 5. The
proof of the main local equivalence result is deferred to Section 6.

\section{Preliminaries}

\subsection{Diffusion processes}

We assume that a continuous record $X^T=\{X_t,\;0\leq t\leq T\}$
of a $d$-dimensional diffusion
process $X$ is observed up to time instant $T$. This diffusion
process is supposed to be given as a solution of the
stochastic differential equation
\begin{equation}\label{diff}
dX_t= b(X_t)\,dt+dW_t,\qquad X_0=\xi,\qquad t\in [0,T],
\end{equation}
where $b:\RR^d\to \RR^d$, $W=(W_t,\,t\geq 0)$ is a $d$-dimensional
Brownian motion and $\xi$ is a random vector independent of $W$. We
denote by $b_i:\RR^d\to\RR$, $i=1,\ldots,d$, the components of the
vector valued function $b$. In what follows, we assume that the
drift is of the form $b=-\nabla V$, where  $V\in C^2(\RR^d)$ is
referred to as potential. This restriction permits to use strong
analytical results for the Markov semigroup of the diffusion on the
$L^2$-space generated by the invariant measure.

For positive constants $M_1$ and $M_2$, we define $\Sigma(M_1,M_2)$
as the set of all functions $b=-\nabla V:\RR^d\to\RR^d$ satisfying
for any $x,y\in\RR^d$
\begin{eqnarray}
|b(x)|&\leq M_1(1+|x|),\label{M}\\
(b(x)-b(y))^T(x-y)&\leq -M_2|x-y|^2,\label{K}
\end{eqnarray}
where $|\cdot|$ denotes the Euclidian norm in $\RR^d$. Any such
function $b$ is locally Lipschitz-continuous. Therefore equation
(\ref{diff}) has a unique strong solution, which is a homogeneous
continuous Markov process, cf. Rogers and Williams \cite{RogWil},
Thm. 12.1. Set $C_b=\int_{\RR^d}e^{-2V(u)}\,du$ and
$$
\mu_b(x)=C_b^{-1}e^{-2V(x)},\qquad x\in\RR^d.
$$
Under condition (\ref{K}) we have $C_b<\infty$ and the process $X$
is ergodic with unique invariant probability measure
(Bhattacharya~\cite[Thm.~3.5]{Bhat}). Moreover, the invariant
probability measure of $X$ is absolutely continuous with respect to
the Lebesgue measure and its density is $\mu_b$. From now on, we
assume that the initial value $\xi$ in (\ref{diff}) follows the
invariant law such that the process $X$ is strictly stationary. We
denote by $\Pb_b^{T}$ the law of this process induced on the
canonical space $\big(C([0,T];\RR^d),\mathcal
B_{C([0,T];\RR^d)}\big)$ and by $\Ex_b$ the expectation operator
with respect to this law.
We write $\mu_b(f):=\Ex_b[f(X_0)]=\int f\mu_b$. Let $P_{b,t}$ be the
transition semigroup of this process on $L^2(\mu_b)$, that is
$$
P_{b,t}f(x)=\Ex_b[f(X_t)|X_0=x],\; f\in
L^2(\mu_b)=\Big\{f:\RR^d\to\RR:\int\!\! |f|^2\mu_b<\infty\Big\}.
$$
The transition density is denoted by $p_{b,t}$: $P_{b,t}f(x)=\int
f(y)p_{b,t}(x,y)\,dy$.

\subsection{Estimators of drift and invariant density}

\subsubsection*{\bf Some notation.}
We write $A(p)\lesssim B(p)$ when $A(p)$ is bounded by a constant
multiple of $B(p)$ uniformly over the parameter values $p$, that is
$A(p)={\cal O}(B(p))$ using the Landau symbol. Similarly,
$A(p)\thicksim B(p)$ means that $A(p)\lesssim B(p)$ as well as
$B(p)\lesssim A(p)$. We denote by $|A|$ the Lebesgue measure and by
$\diam(A)$ the diameter of a Borel set $A\subset\RR^d$.

For any multi-index $\alpha\in\NN^d$ and $x\in\RR^d$ we set
$|\alpha|=\alpha_1+\ldots+\alpha_d$ and
$x^\alpha=x_1^{\alpha_1}\cdot\ldots\cdot x_d^{\alpha_d}$. Let us
introduce the H\"older class
$$
\Hbl=\bigg\{f\in C^{\lfloor\beta\rfloor}(\RR^d;\RR)\,:\;
\begin{matrix}
|D^{\alpha}f(x)-D^{\alpha}f(y)|\leq
L|x-y|^{\beta-\lfloor\beta\rfloor}\\[3pt]
\hbox{for any $\alpha$ such that $|\alpha|=\lfloor \beta\rfloor$}
\end{matrix}
\bigg\}
$$
where $\lfloor\beta\rfloor$ is the largest integer {\it strictly}
smaller than $\beta$ and $D^\alpha f:=\frac{\partial^{|\alpha|}
f}{\partial x_1^{\alpha_1}\ldots\partial x_d^{\alpha_d}}$.

\subsubsection*{\bf The construction.}
Let us assume that the potential $V$ lies in ${\cal H}(\beta+1,L)$
for some $\beta,\,L>0$, which implies $b_i\in {\cal H}(\beta,L)$.
Furthermore, if for some constant $C_1>0$ we have
\begin{equation}\label{deriv}
\max_{i=1,\ldots,d}\max_{\alpha:|\alpha|\leq \lfloor\beta\rfloor}|D^\alpha b_i(0)|\leq C_1
\end{equation}
then the function $\mu_b$ is H\"older continuous of order $\beta+1$ in any bounded set $A\subset\RR^d$, that is
$$
|D^\alpha \mu_b(x)-D^\alpha \mu_b(y)|\leq L_\mu|x-y|^{\beta-\lfloor\beta\rfloor},\qquad \forall
\alpha\in\NN^d:|\alpha|=\lfloor\beta\rfloor+1
$$
for all $x,y\in A$ and for some constant $L_\mu$. We denote by $\widetilde{\mathcal H}(\beta,L,C_1)$
the set of all functions $b$ such that $b_i\in\mathcal H(\beta,L)$ and (\ref{deriv}) is fulfilled.

A natural kernel estimator for the invariant density based on the
observation $X^T$ is given by
\begin{equation}\label{muhat}
\hat\mu_{h,T}(x)=\frac{1}{T}\int_0^T
K_h(x-X_t)\,dt,\quad x\in\RR.
\end{equation}
Here, $K_h(x)=h^{-d}K(h^{-1}x)$ and $K:\RR^d\to\RR$ is a smooth
kernel function of compact support, satisfying $\int K(x)\,dx=1$ and
$\int K(x)x^\alpha\,dx=0$ whenever $1\le
|\alpha|\le\lfloor\beta\rfloor+1$. The usual bias-variance
decomposition and approximation inequality yield
(Efromovich~\cite{EFR}, $\S$ 8.9)
\begin{equation}\label{muhatest}
\Ex_b\big[|\hat{\mu}_{h,T}(x)-\mu_b(x)|^2\big]\lesssim
h^{2(\beta+1)}+T^{-2}\Var\Big[\int_0^T
K_h(x-X_t)\,dt\Big].
\end{equation}
By analogy with the model of regression with random design, a
reasonable estimator of $b$ is obtained by
setting
\begin{equation}\label{bihat}
\hat
b_{h,T}(x)=\frac{\int_0^T K_{h}(x-X_t)\,
dX_t}{T\max(\hat\mu_{h,T}(x),\mu_\ast(x))},\quad
x\in\RR,
\end{equation}
where $\mu_\ast(x)>0$ is some a priori lower bound on $\mu_b(x)$,
see Remark~\ref{lowmu} below. A similar risk analysis gives for
$i=1,\ldots,d$:
\begin{eqnarray}\label{bihatest}
\Ex_b\big[|\hat{b}_{i,h,T}(x)-b_i(x)|^2\big]&\lesssim&
h^{2\beta}+\frac1{Th^d}+\frac1{T^2}\Var\Big[\int_0^T
K_{h}(x-X_t)b_i(X_t)\,dt\Big]\nonumber\\
&&+\Ex_b\big[|\hat{\mu}_{h,T}(x)-\mu_b(x)|^2\big].
\end{eqnarray}

\subsubsection*{\bf Asymptotic results.}
In order to determine the asymptotic behaviour for $T\to\infty$, we
study the variance of general additive functionals of $X$ in $d$
dimensions. To do so, we assume that the semigroup $P_{b,t}$ enjoys
the following properties.

\begin{assumption}[spectral gap inequality]\label{A1}
There exists a $\rho>0$ such that for any $f\in L^2(\mu_b)$ and for
any $t>0$
$$
\|P_{b,t}f-\mu_b(f)\|_{\mu_b}\leq e^{-t\rho}\|f\|_{\mu_b}.
$$
\end{assumption}

\begin{assumption}\label{A2}
There is a $C_0>0$ such that for any $t>0$ and for any pair of
points $x,y\in\RR^d$, satisfying $|x-y|^2<t$, we have
$$
p_{b,t}(x,y)\leq C_0(t^{-d/2}+t^{3d/2}).
$$
\end{assumption}

\begin{remark}
Due to Remark 4.14 in Chen and Wang~\cite{ChenWang}
Assumption \ref{A1} is fulfilled with $\rho=M_2$, whenever (\ref{K})
holds.
\end{remark}

\begin{remark}
If $b$ fulfills (\ref{M}), then Assumption \ref{A2} can be deduced
from Qian and Zheng \cite[Thm.~3.2]{QianZheng}. Indeed, taking in
that inequality $q=1+t$ and bounding the terms $\zeta_q$ and
$\rho_q$ respectively by $Cq^{3/2}$ and $Cq$, we get the desired
inequality. If moreover $b$ is bounded, Assumption \ref{A2} is
satisfied for every $(x,y)\in\RR^d$ and without the term $t^{3d/2}$
at the right-hand side, cf.\ Qian {\it et al.} \cite[inequality
(5)]{QianRusZheng}.
\end{remark}

\begin{proposition}\label{adfunc}
Let $r$ be a positive number and $f:\RR^d\to\RR$ be a bounded,
measurable function with support $\cal S$ satisfying $\diam(|{\cal
S}|)^d<r^d|{\cal S}|$ and $|{\cal S}|<1$. Under Assumptions \ref{A1}
and \ref{A2} there exists a constant $C$ depending only on $r$,
$d\ge 2$ and on $C_0$ and $\rho$ from Assumptions \ref{A1} and
\ref{A2} such that
$$
\Var_b\bigg(\int_0^T f(X_t)\,dt\bigg)\leq C T \|f\|_\infty^2
\mu_b({\cal S})|{\cal S}|\psi_d^2(|{\cal S}|),
$$
where $\|f\|_\infty=\sup_{x\in\RR^d}|f(x)|$ and
$$
\psi_d(x)=\begin{cases} \max(1,(\log(1/x))^2),& d=2,\\
x^{1/d-1/2}, &d\ge 3.\end{cases}
$$
\end{proposition}

\begin{proof}
Set $f_c=f-\mu_b(f)$. Symmetry and stationarity yield
\begin{align*}
\Var_b\bigg(\int_0^T
f(X_t)\,dt\bigg)&=2\int_0^T\int_0^s\Ex_b\big[f_c(X_t)f_c(X_s)\big]dt\,ds\\
&=2\int_0^T\int_0^s\Ex_b\big[f_c(X_0)f_c(X_{s-t})\big]dt\,ds\\
&=2\int_0^T(T-u)\Ex_b\big[f_c(X_0)f_c(X_{u})\big]\,du\\
&\le 2T\int_0^T \big\langle f_c,P_{b,u}f_c\big\rangle_{\mu_b}\,du.
\end{align*}
Let $0<\delta<D\le T$ where the specific choice of $\delta,\,D$ is
given later. Then
\begin{equation}\label{integralgap}
\int_{[0,\delta]\cup [D,T]} \big\langle
f_c,P_{b,u}f_c\big\rangle_{\mu_b}\,du \le (\delta+\rho^{-1}e^{-\rho
D})\|f\|_{\mu_b}^2 \lesssim (\delta+e^{-\rho D})\mu_b({\cal
S})\|f\|_\infty^2
\end{equation}
follows from $\|P_{b,u}f_c\|_{\mu_b}\le e^{-\rho u}\|f\|_{\mu_b}$
given by Assumption \ref{A1}. For moderate values $u\in [\delta,D]$
we use
\begin{align}
\langle f_c, P_{b,u}f_c\rangle_{\mu_b} \le\langle f,
P_{b,u}f\rangle_{\mu_b} \leq \int |f(x)|\Big(\int
p_{b,u}(x,y)\,|f(y)|\,dy\Big)\,\mu_b(x)\,dx \nonumber.
\end{align}
For $\delta>\diam({\cal S})^2$ we infer from Assumption \ref{A2}
\begin{align}\label{deltaDelta}
\langle f, P_{b,u}f\rangle_{\mu_b}&\leq
C(u^{-d/2}+u^{3d/2})\mu_b(|f|)\int |f(y)|\,dy\quad \forall\,u\ge
\delta.
\end{align}
Combining \eqref{integralgap} and \eqref{deltaDelta} and
assuming $\diam({\cal S})<\delta^{1/2}$, for $d>2$ we find
\[ \int_0^T
\big\langle f_c,P_{b,u}f_c\big\rangle_{\mu_b}\,du \lesssim
\Big(\delta+e^{-\rho D}+\delta^{1-d/2}|{\cal S}|+D^{1+3d/2}|{\cal
S}|\Big)\mu_b({\cal S})\|f\|_\infty^2.
\]
Balancing the terms, we choose $D=\max(-\rho^{-1}\log(|{\cal S}|),r^2)$ and
$\delta=r^2|{\cal S}|^{2/d}$. This gives the asserted estimate
because we had assumed $\diam({\cal S})<r|{\cal S}|^{1/d}$. The case
$d=2$ can be treated similarly.
\qed

\end{proof}

\begin{remark}
In the case $d=1$ the bound holds with $\psi_1(x)=1$, cf.
Proposition 5.1 in Dalalyan and Rei{\ss} \cite{DalReiss}.
\end{remark}
\vspace{-10pt}
\begin{remark}
The dimensional effect is due to the singular behaviour of
$p_{b,t}(x,y)$ for $t\to 0$. However, if the term $t^{3d/2}$ is
absent in Assumption \ref{A2}, then in the definition of $\psi_2$
the term $(\log(1/|{\cal S}|))^2$ can be replaced by $(\log(1/|{\cal
S}|))^{1/2}$. This is the case  when the drift is bounded.
\end{remark}

\begin{corollary}\label{estrates}
If $b\in
\widetilde{\cal H}(\beta,L,C_1)\cap \Sigma(M_1,M_2)$, the estimators given in \eqref{muhat} and
\eqref{bihat} satisfy for $h$ sufficiently small the following risk
estimates:
\begin{align*}
\Ex_b\big[(\hat{\mu}_{h,T}(x)-\mu_b(x))^2\big]&\lesssim
h^{2(\beta+1)}+T^{-1}\psi_d^2(h^d),\\
\Ex_b\big[|\hat{b}_{h,T}(x)-b(x)|^2\big]&\lesssim
h^{2\beta}+T^{-1}h^{-d}+h^{2(\beta+1)}+T^{-1}\psi_d^2(h^d).
\end{align*}
The rate-optimal choice $h=h(T)\thicksim T^{-1/(2\beta+d)}$ yields
the rates
\begin{align*}
\Ex_b\big[(\hat{\mu}_{h(T),T}(x)-\mu_b(x))^2\big]^{1/2}&\lesssim
\begin{cases} T^{-1/2}(\log T)^2,&d=2,\\
T^{-(\beta+1)/(2\beta+d)},& d\ge 3,\end{cases}\\
\Ex_b\big[|\hat{b}_{h(T),T}(x)-b(x)|^2\big]^{1/2}&\lesssim
T^{-\beta/(2\beta+d)}.
\end{align*}
\end{corollary}

\begin{proof}
The risk bound for $\hat{\mu}_{h,T}$ follows from
$|\supp(K_h)|\thicksim h^d$, $\|\mu_b\|_\infty\lesssim 1$ and an
application of Proposition \ref{adfunc} to the bias-variance
decomposition \eqref{muhatest} for any $h$ sufficiently small. In
the same way, we obtain the estimate for each $\hat{b}_{i,T,h}$ and
the rates follow by simple substitution.\qed
\end{proof}

\begin{remark}
The convergence rates for the risk of $\hat{\mu}$ are to be compared
with the one-dimensional case, where the parametric rate $T^{-1/2}$
is obtained, and with standard multivariate density estimation,
where the corresponding rate is $n^{-\beta/(2\beta+d)}$ for $n$
observations, which is considerably larger. In contrast, the rate
for $\hat{b}$ corresponds exactly to the classical rate
$n^{-\beta/(2\beta+d)}$ in regression or density estimation.
\end{remark}

\begin{remark}\label{lowmu}
Using conditions (\ref{M}), (\ref{K}) and the equality $
V(x)=V(0)-\int_0^1 b(tx)^Tx\,dt, $ we find
$$
-M_1|x|+\frac12 M_2|x|^2\leq V(x)-V(0)\leq
\frac12 M_1|x|^2+M_1|x|.
$$
Therefore, we can take $\mu_*(x)=e^{-M_1|x|^2-2M_1|x|}/\int e^{2M_1|y|-M_2|y|^2}dy$
as an a priori lower bound for $\mu_b(x)$.
Moreover, due to assumption (\ref{deriv})
the function $\mu_b$ is H\"older continuous in $A_\delta=\{x\in\RR^d: \inf_{y\in A} |x-y|\leq \delta\}$
for any $\delta>0$ and for any bounded set $A\subset\RR^d$. Therefore we do not need to modify the
kernel estimators at the boundaries of $A$ and the inequalities of Corollary~\ref{estrates} hold
uniformly in $b$ and in $x\in A$.
\end{remark}

\begin{remark}\label{supnorm}
Corollary~\ref{estrates} describes the rates of convergence of
estimators for the local risk, that is for a pointwise loss
function. To attain the local neighbourhood defined in the next
section, the risk given by the sup-norm loss must be studied. In the
classical problems of nonparametric estimation, the rates of
convergence for the sup-norm loss on a compact set coincide up to a
logarithmic factor with the local rates of convergence (Korostelev
and Nussbaum~\cite{NK}, Gin\'e, Koltchinskii and Zinn~\cite{GKZ}).
The extension from the pointwise to the uniform loss result is
usually fairly standard, but more involved and lies out of the scope
of this paper.
\end{remark}

\section{Equivalence with the Gaussian shift model}

\subsection{Statement of the result}

Let $\Sigma_\beta(L,M_1,M_2)$ be the set of functions
$b\in\Sigma(M_1,M_2)$ such that all $d$ components $b_i$ of $b$ are
in $\Hbl$. We fix a function $b^{\circ}\in\Sigma_\beta(L,M_1,M_2)$.
Our main result establishes a local asymptotic equivalence between
diffusion and Gaussian shift models in the local setting, that is
when the parameter set is a shrinking neighbourhood of $b^{\circ}$.
${\mathcal B}_E$ always denotes the Borel $\sigma$-algebra of a
topological space $E$.

\begin{definition}[\bf diffusion experiment]
Suppose $\Sigma\subset \Sigma(M_1,M_2)$ for some $M_1,M_2>0$. For
any $T>0$ let $\EE(\Sigma,T)$ be the statistical experiment of
observing the diffusion defined by (\ref{diff}) with $b\in\Sigma$,
that is
$$
\EE(\Sigma,T)=\big(C([0,T];\RR^d), \mathcal B_{C([0,T];\RR^d)},
(\Pb_b^T)_{b\in\Sigma}\big).
$$
\end{definition}

For any function $b\in
L^2(\mu_{b^\circ};\RR^d)=\{f:\RR^d\to\RR^d\,:\,\int
|f|^2\mu_{b^\circ}<\infty\}$ we denote by $\Qb_{b,T}$ the Gaussian
measure on $(C(\RR^d;\RR^d), \mathcal B_{C(\RR^d;\RR^d)})$ induced
by the $d$-dimensional process $Z$ satisfying
\begin{equation}\label{Gauss}
dZ(x)=b(x)\sqrt{\mu_{b^{\circ}}(x)}\,dx+T^{-1/2}\,dB(x),\qquad Z({
\bf 0})={\bf 0},\qquad x\in\RR^d,
\end{equation}
where $B(x)=(B_1(x),\ldots,B_d(x))$ and $B_1(x),\ldots,B_d(x)$ are
independent $d$-variate Brownian sheets, that is zero mean Gaussian
processes with $\cov(B_i(x),B_i(y))=|R_x\cap R_y|$ where
$R_x=\{u\in\RR^d: u_i\in[0,x_i]\}$.

\begin{definition}[\bf Gaussian shift experiment]
For  $\Sigma\subset L^2(\mu_{b^\circ};\RR^d)$ and $T>0$ let
$\FF(\Sigma,T)$ be the Gaussian shift experiment (\ref{Gauss}) with
$b\in\Sigma$, that is
$$
\FF(\Sigma,T)=\big(C(\RR^d;\RR^d), \mathcal B_{C(\RR^d;\RR^d)},
(\Qb_{b,T})_{b\in\Sigma}\big).
$$
\end{definition}

For any positive numbers $\varepsilon$, $\eta$ and for any hypercube
$A\subset\RR^d$, we define the local neighbourhood of $b^\circ$
$$
\Sigma(b^{\circ},\varepsilon,\eta,A)=\bigg\{b\in\Sigma_\beta(L,M_1,M_2):
\begin{matrix}
|b(x)-b^{\circ}(x)|\leq \varepsilon \1_{A}(x),\ x\in\RR^d,\\
|\mu_b(x)-\mu_{b^{\circ}}(x)|\leq \eta\mu_{b^{\circ}}(x),\ x\in A
\end{matrix}
\bigg\},
$$
where $\1_A$ is the indicator function of the set $A$. We state the
main local equivalence result, which will be proved in
Section~\ref{Appendix}. The main ideas of the proof are explained in
the next subsection. For the exact definition of statistical
equivalence and the Le Cam distance $\Delta$ we refer to Le Cam and
Yang~\cite{LCY}.

\begin{theorem}\label{th1}
If $\varepsilon_T$ and $\eta_T$ satisfy the conditions
$$
\lim_{T\to\infty} T^{-\beta}\varepsilon_T^{2-d}= \lim_{T\to\infty}
T^{\frac14+\frac{d-2}{8\beta}}\varepsilon_T(\log(T\varepsilon_T^{-1}))^{\1(d=2)}= \lim_{T\to\infty}
T\eta_T\varepsilon^2_T=0,
$$
then the diffusion model \eqref{diff} is asymptotically equivalent
to the Gaussian shift model (\ref{Gauss}) over the parameter set
$\Sigma_{0,T}=\Sigma(b^{\circ},\varepsilon_T,\eta_T,A)$, that is
$$
\lim_{T\to\infty} \sup_{b^{\circ}\in\Sigma_\beta(L,M_1,M_2)}
\Delta\big(\EE(\Sigma_{0,T},T),\FF(\Sigma_{0,T},T)\big)=0.
$$
\end{theorem}

Let us see for which H\"older regularity $\beta$ on the drift an
estimator can attain the local neighbourhood, that is
$|\hat{b}_{h(T),T}(x)-b(x)|\le\eps_T$ and
$|\hat\mu_{h(T),T}(x)-\mu(x)|\le\eta_T$ hold with a probability
tending to one (cf. Nussbaum \cite{Nussbaum} for this concept). By
the rates obtained in Corollary \ref{estrates}, with a glance at
Remark~\ref{supnorm} and the condition in Theorem \ref{th1}, this is
the case if
\begin{align*}
-\beta-(2-d)\beta/(2\beta+d)&<0,\\
1/4+(d-2)/(8\beta)-\beta/(2\beta+d)&<0,\\
1-(\beta+1)/(2\beta+d)-2\beta/(2\beta+d)&<0.
\end{align*}
It turns out that the second condition is most binding and all three
conditions are satisfied if $\beta>(d-1+\sqrt{2(d-1)^2-1})/2$. The
critical regularity thus grows like $(1/2+1/\sqrt{2})d$ for
$d\to\infty$. In dimension $2$ we obtain the condition $\beta>1$ as
in the result by Carter \cite{Carter} for Gaussian regression.
Whether for H\"older classes of smaller regularity asymptotic
equivalence fails, remains a challenging open problem.

\subsection{Method of proof}\label{Outline}

The general idea of the proof of Theorem~\ref{th1} consists in
discretising (in space) the diffusion process such that the design
regularisation technique we introduced in~\cite{DalReiss} is
applicable in spirit, even though the local time does not exist.

\subsubsection*{\bf Space discretisation.}
For any multi-index $\alpha\in\NN^d$ set $\alpha!=\alpha_1
!\cdot\ldots\cdot\alpha_d!$. Let us denote by
$\{v_i\}_{i=1,\ldots,K}$ the elements of the set
$\{v\in\RR[x]\,:\,v(x)=x^\alpha\hbox{ with }|\alpha|\leq
\lfloor\beta\rfloor\}$ somehow enumerated:
$v_i(x)=x_1^{\alpha_1(i)}\cdot\ldots\cdot
x_d^{\alpha_d(i)}=x^{\alpha(i)}$. We assume that $A=[-a,a[^d$ is a
hypercube and for some $h>0$ with $a/h\in\NN$ we denote by
$\{a_m\}_{m=1,\ldots,M}$ the elements of the grid $(h\ZZ^d)\cap A$.
We introduce the subcubes
$\bC_m=\prod_{j=1}^d[a_{mj},a_{mj}+h[\subset A$, $m=1,\ldots,M$,
where $a_{mj}$ is the $j$th coordinate of $a_m$. Let us define
\begin{equation}
\mv(x)=\begin{pmatrix} v_1(x)/\alpha(1)!\\
\vdots\\ v_K(x)/\alpha(K)!
\end{pmatrix},\label{mbmv}
\end{equation}
which gives rise to the definition $\bar b$ of the Taylor
approximation for $b$
$$
\bar b(x)=\sum_{i=1}^K D^{\alpha(i)}b(a_m)\mv_i(x-a_m) \text{ for }
x\in \bC_m,\; m=1,\ldots,M
$$
and $\bar b(x)=b^{\circ}(x)$ for $x\in\RR^d\setminus A$
($D^{\alpha(i)}$ is applied coordinate-wise). Using this notation,
the Taylor formula can be written as
\begin{equation}\label{Taylor}
b(x)=\bar b(x)+\sum_{i:|\alpha(i)|=\lfloor
\beta\rfloor}\Big(D^{\alpha(i)}
b(\zeta)-D^{\alpha(i)}b(a_m)\Big)\;\frac{v_i(x-a_m)}{\alpha(i)!},\quad
x\in \bC_m,
\end{equation}
where $\zeta\in\RR^d$ satisfies $|\zeta-a_m|\leq |x-a_m|$. This
implies that for $V\in \mathcal H(\beta+1,L)$, the estimate
$|b(x)-\bar b(x)|\lesssim h^\beta$ holds.
We write
\[ \th(x)=b(x)-b^{\circ}(x),\quad \bar\th(x)=\bar b(x)-\bar b^\circ(x)
\text{ and }\bth_j(x)=\begin{pmatrix}D^{\alpha(1)} \th_j(x)\\ \vdots\\
D^{\alpha(K)} \th_j(x)\end{pmatrix}
\]
for $j=1,\ldots,d$ and we shall use equivalently $\theta$ and $b$
for referring to the parameter in the local neighbourhood. The
log-likelihood of the experiment defined via $\Pb_{\bar b}^T$ is
given by (see Liptser and Shiryaev~\cite[p.\ 271, (7.62)]{LipShir})
\begin{align}\label{loglikbarb}
\log\frac{d\Pb_{\bar b}^T}{d\Pb_{\bar b^{\circ}}^T}(X^T)=&
\sum_{m=1}^M\sum_{j=1}^d \Big[\bth_j(a_m)^T\hat\eta_{mj}(T)
-\frac12\; \bth_j(a_m)^T\hat{\mathcal J}_m(T)\,\bth_j(a_m)\Big],
\end{align}
where
\begin{align}
\hat\eta_{mj}(T)&=\int_0^T
\1_{\bC_m}(X_t)\mv(X_t-a_m)\,dW_{t,j}\quad\in\RR^K,\nonumber\\
\hat{\mathcal
J}_m(T)&=\int_0^T\1_{\bC_m}(X_t)\mv(X_t-a_m)\mv(X_t-a_m)^T\,dt\quad\in\RR^{K\times
K},\label{hatsigma}
\end{align}
and $W_{t,j}$ denotes the $j$th component of $W_t\in\RR^d$.

\subsubsection*{\bf Design modification.}
Due to the ergodicity of $X$ the law of the
log-likelihood \eqref{loglikbarb} will for large $T$ be well
approximated by
\begin{equation}\label{log-like}
\sum_{m=1}^M\sum_{j=1}^d\Big(\sqrt{T}\,
\bth_j(a_m)^{T}\eta_{mj}-\frac{T}2\; \bth_j(a_m)^T\mathcal
J_m\bth_j(a_m)\Big)
\end{equation}
where $\eta_{mj}\sim {\cal N}(0,\mathcal J_m)$ i.i.d. and
\begin{equation}\label{sigma}
\mathcal
J_m=\int_{\bC_m}\mv(x-a_m)\mv(x-a_m)^T\mu_{b^{\circ}}(x)\,dx.
\end{equation}

 Since
\begin{eqnarray}
\bth_j(a_m)^T\mathcal J_m\bth_j(a_m)
&=&\int_{\bC_m}(\bar b_j(x)-\bar
b^\circ_j(x))^2\mu_{b^{\circ}}(x)\,dx,\label{theta}
\end{eqnarray}
the process (\ref{log-like}) (indexed by $\bth$) has exactly the
same law as the log-likelihood of the Gaussian shift
$$
dZ(x)=\bar b(x)\sqrt{\mu_{b^{\circ}}(x)}\,dx+T^{-1/2}dB(x),\qquad
Z({\bf 0})={\bf 0},\qquad x\in \RR^d.
$$
Under suitable assumptions on the smoothness of $b$, this
last experiment is asymptotically equivalent to (\ref{Gauss}).


It remains to construct the random variables $(\eta_{mj})$ on some
enlargement of the probability space $(C([0,T];\RR^d), \mathcal
B_{C([0,T];\RR^d)},\Pb_b^T)$ such that $T^{-1/2}\hat\eta_{mj}(T)$
and $\eta_{mj}$ are close as random variables. We define the
stopping time
\begin{equation}\label{taum}
\tau_m =\inf\big\{t\in[0,T]:\|\mathcal J_m^{-1/2}\hat{\mathcal
J}_m(t)\mathcal J_m^{-1/2}\|\geq T\big\}\wedge T,
\end{equation}
where the norm of a matrix $A$ is given by
$\|A\|=\sup_{x}(|Ax|/|x|)$.

Let $\beps=(\varepsilon_{mj})_{m,j}$ be a family of independent
standard normal random vectors in $\RR^K$, defined on an enlarged
probability space such that $\beps$ and $X$ are independent.
We set
$$
\eta_{mj}=\frac1{\sqrt T}\;\hat\eta_{mj}(\tau_m) +(\mathcal
J_m-T^{-1}\hat{\mathcal J}_m(\tau_m))^{1/2}\varepsilon_{mj}.
$$
By definition of $\tau_m$ the matrix $\mathcal
J_m-T^{-1}\hat{\mathcal J}_m(\tau_m)$ is nonnegative definite and
its square root is well defined.

\begin{proposition}\label{prop2}
Under the probability measure $\Pb_{b^\circ}^T$ the random vectors
$(\eta_{mj})_{m,j}\subset\RR^K$ are independent and each $\eta_{mj}$
is centred Gaussian with covariance matrix $\mathcal J_m$.
\end{proposition}

\begin{proof}
It suffices to show that for any sequence
$(\lambda_{mj})_{m,j}\subset\RR^K$ we have
$$
\Ex\bigg[\exp\bigg\{\sum_{m,j}\lambda_{mj}^T\eta_{mj}\bigg\}\bigg]
=\exp\bigg\{\frac12\sum_{m,j}\lambda_{mj}^T\mathcal
J_m\lambda_{mj}\bigg\},
$$
where the expectation is taken with respect to $X$ following the law
$\Pb^T_{b^\circ}$ and $\eps_{mj}$ being i.i.d.\ standard normal in
$\RR^K$, independent of $X$.

The verification of this equality is very similar to the proof of
Proposition 2.13 in Dalalyan and Rei\ss~\cite{DalReiss} and is
omitted. \qed
\end{proof}

\section{Equivalence with heteroskedastic Gaussian regression}

The Gaussian experiment in Theorem~\ref{th1} depends on the centre
$b^{\circ}$ of the neighbourhood via $\mu_{b^{\circ}}$. This fact
makes the passage from the local equivalence to a global equivalence
difficult, especially, because even in the one-dimensional case
there is no known variance stabilising transform for (\ref{Gauss}),
cf. Dalalyan and Rei{\ss} \cite{DalReiss}.

We propose here a method of deriving an asymptotically equivalent
experiment independent of $b^{\circ}$ without using the variance
stabilising transform. The idea is to discretise the Gaussian shift
experiment with a ``step of discretisation'' larger than $1/T$.
This method has already been used in~Brown and
Zhao~\cite{BroZhao} for proving the asymptotic equivalence between
regression models with random and deterministic designs.

We adopt the notation from Section \ref{Outline}. In addition, we
introduce the $K\times K$-matrix
$\mV=\int_{[0,1]^d}\mv(x)\mv(x)^T\,dx,$ where $\mv(x)$ is defined by
(\ref{mbmv}). Since $\mV$ is strictly positive and symmetric, the
matrix $\mV^{-1/2}$ is well defined.

\begin{definition}[\bf heteroskedastic Gaussian regression]
Let $\Sigma$ be a subset of $C^{\lfloor\beta\rfloor}(\RR^d;\RR^d)$.
For any $T,h>0$ we define $\GG(\Sigma,h,T)$ as
the experiment of observing
\begin{equation}\label{reg1}
Y_{im}=\begin{pmatrix}h^{|\alpha(1)|}D^{\alpha(1)} b_i\\ \vdots\\
h^{|\alpha(K)|}D^{\alpha(K)} b_i\end{pmatrix}(a_m)
+{\mV}^{-1/2}\frac{\xi_{im}}{\sqrt{Th^d\mu_b(a_m)}}
\end{equation}
for $i=1,\ldots,d,\ m=1,\ldots,M$, where $(\xi_{im})_{i,m}$ is a
family of independent standard Gaussian random vectors in $\RR^K$
and $b\in\Sigma$.
\end{definition}

Note that the observations in this experiment are chosen from
$\RR^{KMd}$ according to a Gaussian measure. Both the mean and the
variance of this measure depend on the parameter $b$ such that the
experiment is heteroskedastic.

\begin{theorem}\label{th2} If the assumptions of Theorem~\ref{th1}
are fulfilled and $h=h_T$ satisfies
$$
\lim_{T\to\infty} Th^{2\beta}_T=\lim_{T\to\infty}
Th^{2}_T\varepsilon_T^2=\lim_{T\to\infty} \eta^2_Th^{-d}_T=0,
$$
then the diffusion experiments  and the heteroskedastic Gaussian
regression experiments are asymptotically equivalent, that is
$$
\lim_{T\to\infty} \sup_{b^{\circ}\in\Sigma_\beta(L,M_1,M_2)}
\Delta\big(\EE(\Sigma_{0,T},T),\GG(\Sigma_{0,T},h_T,T)\big)=0.
$$
\end{theorem}

\begin{proof}
Theorem~\ref{th1} yields the asymptotic equivalence of the
experiment $\EE$ with the (translated) Gaussian shift experiment
\[ d\widetilde
Z(x)=(b-b^\circ)(x)\sqrt{\mu_{b^\circ}(x)}\,dx+T^{-1/2}dB(x),\quad
x\in\RR^d.
\]
Let us introduce a new Gaussian shift:
$$
d\widehat Z(x)=\sum_{m=1}^M\Big((\bar b-
b^\circ)(x)\sqrt{\mu_{b^\circ}(a_m)}\Big)\1_{\bC_m}(x)\,dx+T^{-1/2}dB(x),\quad
x\in\RR^d.
$$
Since $|\nabla \mu_b(x)|$ and $|\mu_b(x)|$ are uniformly bounded,
the difference between the drifts of $\widetilde Z$ and $\widehat Z$
can be estimated as follows:
\begin{align*}
\big|(b&-b^\circ)(x)\sqrt{\mu_{b^\circ}(x)}-(\bar b-
b^\circ)(x)\sqrt{\mu_{b^\circ}(a_m)}\big|\\
&\leq\big|(b-\bar b)(x)\sqrt{\mu_{b^\circ}(a_m)}\big|
+\big|(b-b^\circ)(x)\big(\sqrt{\mu_{b^\circ}(x)}-\sqrt{\mu_{b^\circ}(a_m)}\big)\big|\\
&\lesssim h^{\beta}+\varepsilon h\qquad \forall x\in \bC_m.
\end{align*}
Therefore, the Hellinger distance between the measures induced by
$\widetilde Z$ and $\widehat Z$ tends to zero as $T\to\infty$
(Strasser~\cite[Rem. 69.8.(2)]{Str}), provided that
$T\varepsilon^2h^2\to 0$ and $Th^{2\beta}\to 0$. The log-likelihood
of the experiment given by $\widehat Z$ has exactly the same law as
the log-likelihood of the Gaussian regression
\begin{equation}\label{reg1'}
Y_{im}=\begin{pmatrix}h^{|\alpha(1)|}D^{\alpha(1)} b_i\\ \vdots\\
h^{|\alpha(K)|}D^{\alpha(K)} b_i\end{pmatrix}(a_m)
+{\mV}^{-1/2}\frac{\xi_{im}}{\sqrt{Th^d\mu_{{b^{\circ}}}(a_m)}}
\end{equation}
for $i=1,\ldots,d;\ m=1,\ldots,M$, where $(\xi_{im})_{i,m}$ is a
family of independent standard Gaussian random vectors in $\RR^K$
and $b\in\Sigma$. By Lemma~3 from Brown et {\it al.}~\cite{BCLZ} the
square of the Hellinger distance between the measures induced by the
observations (\ref{reg1}) and (\ref{reg1'}), respectively, is up to
a constant bounded by
$\sum_{m=1}^M(\mu_b(a_m)-\mu_{b^{\circ}}(a_m))^2/\mu_{b^{\circ}}(a_m)^2\lesssim
M\eta_T^2$. Because of $Mh^d=|A|$ we infer $M\thicksim h^{-d}$ and
the condition $h_T^{-d}\eta_T^2\to 0$ as $T\to\infty$ implies that
the Hellinger distance tends to zero uniformly in
$b\in\Sigma_{0,T}$. Finally, the desired result follows by bounding
the Le Cam distance between experiments by the supremum of the
Hellinger distance between the corresponding measures, see e.g.
Nussbaum~\cite[Eq. (12)]{Nussbaum}.\qed
\end{proof}

\begin{remark}
The experiment given by~(\ref{reg1}) is more informative than the
experiment generated by the observations $(\be_1^TY_{im})_{i,m}$,
where $\be_1=(1,0,\ldots,0)^T\in\RR^K$. If we enumerate
$\{\alpha(i)\}_i$ so that $\alpha(1)={\bf 0}\in\RR^d$ then
$\widetilde Y_{m}:=(\be_1^TY_{1m},\ldots,\be_1^TY_{dm})^T$ satisfies
$\widetilde Y_{m}= b(a_m)+{\epsilon_{m}}/{\sqrt{Th^d\mu_b(a_m)}}$
with $\epsilon_{m}/\sqrt{(\mV^{-1})_{11}}\sim \mathcal N(0,I_d)$
i.i.d. Therefore the diffusion experiment $\EE(\Sigma_{0,T},T)$ is
asymptotically more informative than the regression experiment:
$$
\widetilde Y_{m}=
b(a_m)+\frac{\epsilon_{m}}{\sqrt{Th^d\mu_b(a_m)}},\quad
m=1,\ldots,M.
$$
\end{remark}

If we choose $h_T=T^{-\alpha}$,
$\varepsilon_T=T^{-\beta/(2\beta+d)}$ and
$\eta_T=T^{-(\beta+1)/(2\beta+d)}$ (in view of
Corollary~\ref{estrates}), the condition of Theorem~\ref{th2} takes
the form
$$
\max\bigg(\frac1{\beta};\;\frac{d}{2\beta+d}\bigg)<2\alpha<\frac{4(\beta+1)}{d(2\beta+d)}.
$$
Such a value $\alpha$ exists if and only if
$$
\beta>\max\bigg(\frac{d^2}4-1;\;\frac{d-2+\sqrt{(d-2)^2+4d^2}}{4}\bigg).
$$
For $d=2$ this inequality reduces to $\beta>1$. For $d\geq 4$ it is
equivalent to $\beta>(d/2)^2-1$. Note also that the
logarithmic factors in $\varepsilon_T$ and $\eta_T$ do not affect
this bound on the minimal regularity.

As mentioned in the introduction, the result of Theorem~\ref{th2} is
new already in the one-dimensional case. When $d=1$, using a
$\sqrt{T}$-consistent estimator of $\mu_b$
(Kutoyants~\cite{Kutoyants}, $\S$ 4.2), the local neighbourhood can
be attained as soon as $\beta>1/2$. Taking $K=1$ and using the
globalisation method developed in~\cite{DalReiss}, we obtain the
global asymptotic equivalence of the diffusion experiment and the
regression
$$
Y_m=b(a_m)+\frac{\epsilon_m}{\sqrt{Th\mu_b(a_m)}},\qquad m=1,\ldots,M,
$$
provided that $h=h_T=T^{-\alpha}$ with $ (2\beta)^{-1}<\alpha<1$ and
the assumptions of \cite[Thm. 3.5]{DalReiss} are fulfilled.

\section{Equivalence mapping}

The result of Theorem~\ref{th1} implies in particular that there
exists a Markov kernel $K$ from $(C([0,T];\RR^d),\mathcal
B_{C([0,T];\RR^d)})$ to $(C(\RR^d;\RR^d),\mathcal
B_{C(\RR^d;\RR^d)})$ such that
$$
\lim_{T\to\infty} \sup_{b\in
\Sigma_{0,T}}\|\Pb_b^TK-\Qb_{b,T}\|_{TV}= 0,
$$
where $\Pb_b^TK(A)=\int_{C([0,T];\RR^d)}K(x,A)\Pb_b^T(dx)$ for $A\in
\mathcal B_{C(\RR^d;\RR^d)}$ and $\|\cdot\|_{TV}$ denotes the total
variation norm. The aim of this section is to construct this Markov
kernel explicitly. The construction is divided into two steps.
First, we give the Markov kernel from the diffusion experiment to a
suitable multivariate Gaussian regression. Then we give the Markov
kernel from the Gaussian regression to the Gaussian shift
experiment. An explicit Markov kernel in the other direction is not
known, but seems also less useful.

Assume that we have a path $X^T$ of the diffusion
process~(\ref{diff}) at our disposal. In what follows we use the
notation introduced in Section~\ref{Outline} with $h$ verifying
(\ref{hrate}) below. For any $i=1,\ldots,d$ we denote by $X_{t,i}$
the $i$th coordinate of $X_t$ and define the randomisation
\begin{eqnarray*}
\Phi^{(1)}_{im}(X^T,\beps)&=&\frac1T\int_0^{\tau_m}\1_{\bC_m}(X_t)
\mv(X_t-a_m)\,(dX_{t,i}-\bar b_i^\circ(X_t)\,dt)\\
&&+ \frac1{\sqrt T}(\mathcal J_m-T^{-1}\hat{\mathcal
J}_m(\tau_m))^{1/2}\varepsilon_{im},\quad m=1,\ldots,M,
\end{eqnarray*}
where $\hat{\mathcal J}_m(t)$, $\mathcal J_m$ and $\tau_m$ are
defined by (\ref{hatsigma}), (\ref{sigma}) and (\ref{taum}) and
$\beps=(\varepsilon_{im})_{i,m}$ is a family of independent (and
independent of $X^T$) standard Gaussian vectors in $\RR^K$. As is
easily checked, the random vector $\mathcal
J_m^{-1}\Phi^{(1)}_{im}(X^T,\tilde\beps)$ with
$\tilde\varepsilon_{im}=(T\mathcal J_m-\hat{\mathcal
J}_m(\tau_m))^{1/2} \bth_i(a_m)+\varepsilon_{im}$ has the same law
as the Gaussian regression
\begin{equation}\label{regr1}
Y_{im}=\bth_i(a_m)+(T\mathcal J_m)^{-1/2}\varepsilon_{im}.
\end{equation}
We prove in Section~\ref{Proof Th1} that the total variation between the laws of $\beps$
and $\tilde\beps$ tends to zero as $T\to\infty$. Consequently, if we
denote by $K^{(1)}(x,\cdot)$ the law of
$\{\mathcal J_m^{-1}\Phi^{(1)}_{im}(x,\beps);\ i=1,\ldots,d;\
m=1,\ldots,M\}$, we obtain a Markov
kernel realising the asymptotic equivalence between the diffusion
\eqref{diff} and the Gaussian regression (\ref{regr1}).

For any $x\in\bC_m$ and for any $i\in\{1,\ldots,d\}$, we define the
randomisation of the regression \eqref{regr1} by
\begin{eqnarray}
&\Phi_{i,x}^{(2)}&(Y,\tilde B)=\int_{R(a_m,x)}\!\!\!\!\!\!\big(\bar
b^\circ_i(u)+\mv(u)^TY_{im}\big)\sqrt{\mu_{b^\circ}(u)}\,du\nonumber\\
&&+\frac1{\sqrt T}\int_{R(a_m,x)}\!\!\!\!\!\!\!\!\!\sqrt{\mu_{b^\circ}(u)}\,d\tilde B_i(u)\nonumber\\
&&-\frac1{\sqrt T}\Big(\int_{R(a_m,x)}\!\!\!\!\!\!\mv(u)^T\mu_{b^\circ}(u)\,du\Big)
\mathcal J_m^{-1}\Big(\int_{\bC_m}\!\!\!\!\!\!\mv(u)\sqrt{\mu_{b^\circ}(u)}\,d\tilde
B_i(u)\Big),\qquad
\label{Tphi2}
\end{eqnarray}
where $R(a_m,x)=\prod_{i=1}^d[a_{mi},x_i[$, $\tilde B=(\tilde
B_1,\ldots,\tilde B_d)$ and $\tilde B_1,\ldots,\tilde B_d$ are
independent $d$-variate Brownian sheets independent of
$(Y_{im})_{i,m}$. Let us show that $\Phi^{(2)}(\mathbf y,\tilde
B)=(\Phi^{(2)}_{i,x}(\mathbf y,\tilde B);\, i\in\{1,\ldots,d\},x\in
A)$ is an equivalence mapping from the Gaussian regression model
(\ref{regr1}) to the Gaussian shift model (\ref{Gauss}).

For any $x\in \bC_m$ and for any $i=1,\ldots,d$ define the multivariate analogue of a Brownian bridge
\begin{align*}
V_i(x)&=\int_{R(a_m,x)} \mv(u)\sqrt{\mu_{b^\circ}(u)}\,d\tilde
B_i(u)\\
&\quad-\Big(\int_{R(a_m,x)}\mv(u)\mv(u)^T\mu_{b^\circ}(u)\,du\Big)\mathcal J_m^{-1}\Big(\int_{\bC_m}
\mv(u)\sqrt{\mu_{b^\circ}(u)}\,d\tilde B_i(u)\Big)
\end{align*}
and set
$$
\tilde
V_i(x)=\Big(\int_{R(a_m,x)}\mv(u)\mv(u)^T\mu_{b^\circ}(u)\,du\Big)\,Y_{im}+T^{-1/2}V_i(x).
$$
The process $\tilde V_i$ takes values in $\RR^K$
and can be rewritten in the form
$\tilde V_i(x)
=\int_{R(a_m,x)}\mv(u)(\bar b_i(u)-\bar
b_i^\circ(u))\mu_{b^\circ}(u)\,du+T^{-1/2}\widehat W_i(x)
$
where
\[ \widehat
W_i(x)=\Big(\int_{R(a_m,x)}\mv(u)\mv(u)^T\mu_{b^\circ}(u)\,du\Big)\mathcal J_m^{-1/2}\varepsilon_{im}+V_i(x).
\]
By construction, the process $\widehat W_i$ is centred Gaussian with
covariance matrix $\Ex[\widehat W_i(x)\widehat W_i(\bar
x)^T]=\int_{R(a_m,x)\cap R(a_m,{\bar
x})}\mv(u)\mv(u)^T\mu_{b^\circ}(u)\,du$. Assuming that
$v_1,\ldots,v_K$ are enumerated in such a way that $v_1(u)\equiv 1$,
one checks that $\widehat
B_i(x)=\int_{R(a_m,x)}\mu_{b^\circ}(u)^{-1/2}d\widehat W_{i,1}(u)$
is a $d$-variate Brownian sheet, where $\widehat W_{i,1}$ is the
first coordinate of $\widehat W_i$. Therefore, the randomisation
\begin{equation}\label{phi21}
\Phi^{(2)}_{i,x}(Y,\tilde B)=
\int_{R(a_m,x)}\bar b^\circ_i(u)\sqrt{\mu_{b^\circ}(u)}\,du+\int_{R(a_m,x)}\mu_{b^\circ}(u)^{-1/2}d\tilde
V_{i,1}(u)
\end{equation}
satisfies
\begin{equation}\label{phi2}
d\Phi^{(2)}_{i,x}=\bar
b_i(x)\,\sqrt{\mu_{b^\circ}(x)}\,dx+T^{-1/2}d\widehat B_i(x),\quad
x\in\bC_m,\quad i=1,\ldots d.
\end{equation}
The total variation between the measures induced by (\ref{phi2}) and
(\ref{Gauss}) is up to a constant bounded by $\sqrt T h^{\beta}$,
which tends to zero because of our choice of $h$ and the assumptions
of Theorem~\ref{th1}.
Moreover, the $d$-variate Brownian sheets $\widehat B_1,\ldots,\widehat B_d$ are independent.
Simple algebra shows that
the two definitions (\ref{phi21}) and (\ref{Tphi2}) coincide.
Hence the law $K^{(2)}(\mathbf y,\cdot)$ of
$\Phi^{(2)}({\mathbf  y},\tilde B)$ provides a Markov
kernel from the Gaussian regression (\ref{regr1}) to the Gaussian
shift (\ref{Gauss}) realising the asymptotic equivalence.

\section{Proof of Theorem \ref{th1}}\label{Appendix}

\subsection{Main part}\label{Proof Th1}

As we have seen in Section \ref{Outline}, the construction of the
Gaussian experiment makes use of an i.i.d.\ family
$\beps=(\varepsilon_{mj})_{m=1,\ldots,M,\;j=1,\ldots,d}$ of standard
Gaussian vectors with values in $\RR^K$. The canonical version of
$\beps$ is defined on the measurable space $\big(\RR^{KMd},\mathcal
B_{\RR^{KMd}}\big)$. We prove the asymptotic equivalence by a
suitable coupling, which consists in constructing probability
measures $\tilde \Pb_b^{T}$ and $\tilde\Qb_b^{T}$ on the product
space
$$
(\mathscr E,\mathcal B_{\mathscr E}):=\big(C([0,T],\RR^d)\times
\RR^{KMd} ,\mathcal B_{C([0,T],\RR^d)}\otimes \mathcal
B_{\RR^{KMd}}\big)
$$
such that\\[-8mm]
\begin{enumerate}
\item[a)] $\EE(\Sigma_{0,T},T)$ is equivalent to
$\tilde\EE(\Sigma_{0,T},T)=\big(\mathscr E, \mathcal B_{\mathscr E},
(\tilde\Pb_b^T)_{b\in\Sigma_{0,T}}\big)$,
\item[b)] $\tilde\EE(\Sigma_{0,T},T)$ and
$\tilde\FF(\Sigma_{0,T},T)=\big(\mathscr E, \mathcal B_{\mathscr E},
(\tilde\Qb_b^T)_{b\in\Sigma_{0,T}}\big)$ are asymptotically
equivalent,
\item[c)] $\FF(\Sigma_{0,T},T)$ is asymptotically equivalent
to $\tilde\FF(\Sigma_{0,T},T)$.
\end{enumerate}

a) Define $\tilde\Pb^T_b$ to be the measure induced by the pair
$(X^T, \beps)$, where $X^T$ is given by (\ref{diff}) and $\beps$ is
a standard Gaussian vector independent of $X^T$, that is
$\tilde\Pb^T_b=\Pb^T_b\otimes \Norm_{KMd}$ with $\Norm_k$ denoting
the standard normal law on $\RR^k$. Then the equivalence
$\EE\sim\tilde\EE$ follows from the equality in law of the
respective likelihood processes, cf.\ Strasser~\cite[Cor.
25.9]{Str}.

b) The measure $\tilde\Qb_b^T$ is defined via
$$
\tilde\Qb_b^T(A\times B)=\int_{A\times B}
e^{f_b(X^T,\beps)}\,\Pb_{b^{\circ}}^T(dX^T)\,\Norm_{KMd}(d\beps)
$$
for $A\in \mathcal B_{C([0,T],\RR^d)}$ and $B\in \mathcal
B_{\RR^{KMd}}$ with
$$
f_b(X^T,\beps)= \sum_{m=1}^M\sum_{j=1}^d
\bigg[\sqrt{T}\bth_j(a_m)^T\eta_{mj}(X^T,\beps) -\frac{T}2\;
\bth_j(a_m)^T\mathcal J_m\,\bth_j(a_m)\bigg]
$$
and
\begin{align*}
\eta_{mj}(X^T,\beps)=&\frac1{\sqrt
T}\int_0^{\tau_m}\1_{\bC_m}(X_t)\mv(X_t-a_m)\,(dX_{t,j}-b_{j}^\circ(X_t)\,dt)\\
&+(\mathcal J_m-T^{-1}\hat{\mathcal
J}_m(\tau_m))^{1/2}\varepsilon_{mj}.
\end{align*}
Because of $f_{b^{\circ}}(X^T,\varepsilon)=0$ these definitions
yield $\tilde\Qb_{b^{\circ}}^T=\tilde\Pb_{b^{\circ}}^T$ and
therefore
$\log\big(\frac{d\tilde\Qb_b^T}{d\tilde\Qb_{b^{\circ}}^T}\;(X^T,\beps)\big)=f_b(X^T,\beps)$.
Proposition~\ref{prop2} combined with the classical formula of the characteristic function of a
Gaussian vector implies that $\tilde\Qb_b^T$ is a probability measure.

To prove the asymptotic equivalence of $\tilde\EE$ and $\tilde\FF$,
it suffices to show that the Kullback-Leibler divergence between the
measures $\tilde\Pb_b^T$ and $\tilde\Qb_b^T$ tends to zero uniformly
in $b\in\Sigma_{0,T}$ (see the proof of Thm. 2.16 in~\cite{DalReiss}).
The Fubini theorem yields
\begin{eqnarray*}
KL(\tilde\Pb_b^T,\tilde\Qb_b^T)&=&\int\log\Big(\frac{d\tilde\Pb_b^T}{d\tilde\Qb_b^T}(X^T,\beps)\Big)
\Pb_b^T(dX^T)\Norm_{KMd}(d\beps)\\
&=&\Ex_b\Big[\log\Big(\frac{d\Pb_b^T}{d\Pb_{b^{\circ}}^T}(X^T)\Big)-\int
f_b(X^T,\beps)\,\Norm_{KMd}(d\beps)\Big].
\end{eqnarray*}
The Girsanov formula (Liptser and Shiryaev \cite{LipShir}) and the
fact that the expectation of the stochastic integral is zero give
\begin{align*}
&\Ex_b\Big[\log\Big(\frac{d\Pb_b^T}{d\Pb_{b^{\circ}}^T}(X^T)\Big)\Big]=
\Ex_b\Big[\log\Big(\frac{\mu_b(X_0)}{\mu_{b^{\circ}}(X_0)}\Big)\Big]+
\frac12\Ex_b\Big[\int_0^T|\th(X_t)|^2 \,dt\Big]\\
&=\Ex_b\Big[\log\Big(\frac{\mu_b(X_0)}{\mu_{b^{\circ}}(X_0)}\Big)\Big]+
\frac{T}2
\,\int_{A}\big|\th(x)-\bar\th(x)\big|^2\,\mu_b(x)\,dx\\
&\quad + \frac{T}2 \,\int_{A}|\bar\th(x)|^2\mu_b(x)\,dx + T
\,\int_{A}\bar\th(x)^T\big(\th(x)-\bar\th(x)\big)\,\mu_b(x)\,dx.
\end{align*}
Similarly, we find
\begin{align*}
&\Ex_b\Big[\int f_b(X^T,\beps)\,\Norm_{KMd}(d\beps)\Big]
=\sum_{m=1}^M\sum_{j=1}^d\Big(-\frac{T}2\;
\bth_{j}(a_m)^T\mathcal J_m\,\bth_j(a_m)\\
&\qquad\qquad\qquad+\Ex_b\Big[\bth_j(a_m)^T
\int_0^{\tau_m}\1_{\bC_m}(X_t)\mv(X_t-a_m)\,\th_{j}(X_t)\,dt\Big]\Big)\\
&\quad =- \frac{T}{2} \,\int_{A}|\bar\th(x)|^2\mu_{b^\circ}(x)\,dx +
\sum_{m=1}^M\Ex_b\Big[\int_0^{\tau_m}\1_{\bC_m}(X_t)|\bar\th(X_t)|^2\,dt\Big]\\
&\qquad\qquad\qquad
+\sum_{m=1}^M\Ex_b\Big[\int_0^{\tau_m}\1_{\bC_m}(X_t)\bar\th(X_t)^T(\th(X_t)
-\bar\th(X_t))\,dt\Big].
\end{align*}
Using for $f(x)=|\bar\th(x)|^2$ and
$f(x)=\bar\th(x)^T\big(\th(x)-\bar\th(x)\big)$ the general identity
\[
T \,\int_{A}f(x)\mu_b(x)\,dx=\sum_{m=1}^M
\Ex_b\Big[\int_0^T\1_{\bC_m}(X_t)\,f(X_t)\,dt\Big],
\]
we obtain $\KL(\tilde\Pb_b^T,\tilde\Qb_b^T)
=\sum_{i=1}^5\calT_i(\th)$ with
\begin{eqnarray*}
\calT_1(\th)&=&\Ex_b\big[\log\mu_b(X_0)-\log\mu_{b^{\circ}}(X_0)\big],\\
\calT_2(\th)&=&\frac{T}2\int_{A}|\bar\th(x)|^2\big(\mu_{b^{\circ}}(x)-\mu_b(x)\big)\,dx,\\
\calT_3(\th)&=&\sum_{m=1}^M\Ex_b\Big[\int_{\tau_m}^T
|\bar\th(X_t)|^2\1_{\bC_m}(X_t)\,dt\Big],\\
\calT_4(\th)&=&\frac{T}2
\,\int_A \big|\th(x)-\bar\th(x)\big|^2\,\mu_b(x)\,dx,\\
\calT_5(\th)&=&\sum_{m=1}^M\Ex_b\Big[\int_{\tau_m}^T\1_{\bC_m}(X_t)\,\bar\th(X_t)^T(\th(X_t)
-\bar\th(X_t))\,dt\Big].
\end{eqnarray*}
The Cauchy-Schwarz inequality implies that $\calT_5(\th)\leq
\calT_3(\th)+\calT_4(\th)$. The explicit form of the invariant
density $\mu_b$ implies that $\sup_\vartheta
\calT_1(\vartheta)\lesssim \varepsilon$. The H\"older assumption
implies that $\sup_x|\bar\vartheta(x)-\vartheta(x)|\lesssim h^\beta$
and we infer
$$
\sup_\th \calT_2(\th)\lesssim T(h^{2\beta}+\varepsilon^2)\eta,\qquad
\sup_\th \calT_4(\th)\lesssim
Th^{2\beta}.
$$
In Section \ref{LemmaT3} below we prove that
\begin{equation}\label{EstT3}
\calT_3(\vartheta)\lesssim(T\eta+\psi_d(h^d)\sqrt T )\,\|\bar
\vartheta\|_\infty^2
\end{equation}
holds if $h=h_T$ tends to zero for $T\to\infty$. Hence, we obtain
$$
KL(\tilde\Pb_b^T,\tilde\Qb_b^T)\lesssim
\varepsilon+Th^{2\beta}+T(\varepsilon^2+h^{2\beta})\eta+\psi_d(h^d)
\sqrt T(\varepsilon^2+h^{2\beta}).
$$
Consequently, the rate-optimal choice of $h$ is
\begin{equation}\label{hrate}
h=h_T=(\varepsilon^4T^{-1})^{1/(4\beta+d-2)},
\end{equation}
provided that $h^{2\beta}=o(\varepsilon^2)$, so that
$$
KL(\tilde\Pb_b^T,\tilde\Qb_b^T)\lesssim
\varepsilon+(\varepsilon^2T^{\frac12+\frac{d-2}{4\beta}})^{4\beta/(4\beta+d-2)}(\log(T\varepsilon^{-1}))^{2\1(d=2)}
+T\varepsilon^2\eta,
$$
given $\varepsilon^{d-2} T^{\beta}\to \infty$. Under the assumptions
of the theorem we thus conclude that $\tilde\EE$ and $\tilde\FF$ are
asymptotically equivalent.

c) It remains to verify that the statistical experiment $\FF$
defined via $\Qb_b^T$ is asymptotically equivalent to the experiment
$\tilde\FF$ defined via $\tilde\Qb_b^T$. We have already seen that
$$
\log\bigg(\frac{d\tilde\Qb_b^T}{d\tilde\Qb_{b^{\circ}}^T}\bigg)=\sum_{m,j}
\bigg[\sqrt{T}\bth_j(a_m)^T\eta_{mj} -\frac{T}2\;
\bth_j(a_m)^T\mathcal J_m\,\bth_j(a_m)\bigg].
$$
Recall that according to Proposition~\ref{prop2} the random vectors
$(\eta_{mj})_{m,j}$ are independent Gaussian with covariance matrix
$\mathcal J_m$. Therefore, the law of the log-likelihood process
$\big(d\tilde\Qb_b^T/d\tilde\Qb_{b^{\circ}}^T\big)_{b\in\Sigma_0}$
coincides with the law of the process $\big({d\tilde\Qb_{\bar
b}^T}/{d\tilde\Qb_{\bar b^{\circ}}^T}\big)_{b\in\Sigma_0}$. This
gives the equivalence of the experiments $\tilde\FF$ and
$\widehat\FF$, where the latter experiment is defined by the
observation
\begin{equation}\label{hatF}
dZ(x)=\bar b(x)\sqrt{\mu_{b^{\circ}}(x)}\,dx+T^{-1/2}\,dB(x),\qquad
Z({\bf 0})={\bf 0},\qquad x\in\RR^d.
\end{equation}
To conclude, we remark that the Kullback-Leibler divergence between
the Gaussian experiments $\FF$ and $\widehat\FF$ is bounded by
$T\int_{\RR^d}(\bar b-b)^2\mu_{b^{\circ}}\leq Th_T^{2\beta}$ and in
view of \eqref{hrate} tends to zero for $T\to \infty$. \qed

\subsection{Evaluation of ${\mathcal T}_3$ }\label{LemmaT3}

We start by sketching how the estimate could be reduced to a purely
analytical problem, using
\begin{align}
&{\mathcal T}_3(\vartheta)\le \norm{\bar{b}-\bar{b}_0}_\infty^2
\sum_m
\Ex_b\Bigl[\int_{\tau_m}^T \1_{\bC_m}(X_t)\,dt\Bigr]\label{T3decomp}\\
&\quad\le \norm{\bar{b}-\bar{b}_0}_\infty^2 \Big(\sup_m
\Ex_b[T-\tau_m]+\sum_m
\Big(\Ex_b\Bigl[\int_{\tau_m}^T(\1_{\bC_m}(X_t)-\Pb_b(\bC_m))\,dt\Bigr]\Big).
\nonumber
\end{align}
If $f$ is a function in the domain of the generator $L_b$ of the
semigroup $(P_{b,t})_{t\ge 0}$ with
$L_bf=\1_{\bC_m}(X_t)-\Pb_b(\bC_m)$, then Dynkin's formula and the
fact that $\1_{\bC_m}(X_t)-\Pb_b(\bC_m)$ is centred yield
\[
\Ex_b\Bigl[\int_{\tau_m}^T(\1_{\bC_m}(X_t)-\Pb_b(\bC_m))\,dt\Bigr] =
\Ex_b[f(X_{\tau_m})]\le \sup_{x\in\bC_m}f(x).
\]
Unfortunately, a suitably tight supremum norm estimate for
$f=L_b^{-1}(\1_{\bC_m}-\Pb_b(\bC_m))$ could not be found in the
literature.

We therefore proceed differently and make use of the mixing
properties of $X$. Fix some $\Delta=\Delta(T)>0$. Since for
$\tau_m>T-\Delta$ the integral over $[\tau_m,T]$ is smaller than the
integral over $[T-\Delta,T]$, we have
\begin{align}\label{intdecomp}
\Ex_b\Bigl[\int_{\tau_m}^T \1_{\{X_t\in {\bC}_m\}}\,dt\Bigr] &\le
\Delta \mu_b(C_m)+ \Ex_b\Bigl[\1_{\{\tau_m\le
T-\Delta\}}\int_{\tau_m}^T \1_{\{X_t\in \bC_m\}}\,dt\Bigr].
\end{align}

\begin{lemma}\label{lem3} Under the assumptions of Proposition
\ref{adfunc} we obtain
$$
\Ex_b\Bigl[\1_{\{\tau_m\le T-\Delta\}}\int_{\tau_m}^{\tau_m+\Delta}
\1_{\{X_t\in \bC_m\}}\,dt\Bigr]\lesssim
\Delta\mu_b(\bC_m)+h^{\frac{d}2}\psi_d(h^d)\sqrt{T \mu_b(\bC_m)}.
$$
\end{lemma}

\begin{proof}
Because of $[\tau_m,\tau_m+\Delta]\subset [(i-1)\Delta,(i+1)\Delta]$
for some $1\le i\le T/\Delta$ we get
$$
\int_{\tau_m}^{\tau_m+\Delta} \1_{\bC_m}(X_s)\,ds\leq
\max_{i=1,\ldots,[T/\Delta]}
\int_{(i-1)\Delta}^{(i+1)\Delta}\1_{\bC_m}(X_s)\,ds.
$$
Set
$U_i=\int_{(i-1)\Delta}^{(i+1)\Delta}\1_{\bC_m}(X_s)\,ds-2\Delta\mu_b(\bC_m)$.
By separating the bias from the stochastic term, we find
$$
\int_{\tau_m}^{\tau_m+\Delta} \1_{\bC_m}(X_s)\,ds\leq 2\Delta
\mu_b(\bC_m)+\max_{i=1,\ldots,[T/\Delta]} |U_i|,
$$
and by the Cauchy-Schwarz inequality
\begin{align*}
\Ex_b\big[\max_{i} |U_i|\big]&\leq
\Big(\sum_{i=1}^{\lfloor T/\Delta\rfloor }
\Ex_b(U_i^2)\Big)^{\frac12}=\lfloor
T/\Delta\rfloor^{1/2}\Var\bigg(\int_{0}^{2\Delta}\1_{\bC_m}(X_s)\,ds\bigg)^{\frac12}.
\end{align*}
We conclude by an application of Proposition~\ref{adfunc}.
\qed\end{proof}

\begin{lemma} If Assumption \ref{A1} is satisfied, then
\begin{align*}
\Ex_b\Bigl[\1_{\{\tau_m\le T-\Delta\}}\int_{\tau_m+\Delta}^T
\1_{\{X_t\in \bC_m\}}\,dt\Bigr]&\leq
\mu_b(\bC_m)\int_0^{T-\Delta}\Pb_b^T(\tau_m\leq
t)\,dt\\
&\quad+Te^{-\Delta\rho}\sqrt{\mu_b(\bC_m)}.
\end{align*}
\end{lemma}
\begin{proof}
We have
\begin{align*}
\Ex_b\Bigl[\1_{\{\tau_m\le T-\Delta\}}\int_{\tau_m+\Delta}^T
&\1_{\{X_t\in \bC_m\}}\,dt\Bigr] =
\Ex_b\Bigl[\int_{\Delta}^T\1_{\{X_t\in \bC_m\}}\,\1_{\{\tau_m\leq
t-\Delta\}}\,dt\Bigr]\\
&=\mu_b(\bC_m)\int_\Delta^T\Pb_b^T(\tau_m\leq
t-\Delta)\,dt\\
&\quad
+\int_{\Delta}^T\Ex_b\bigl[\big(\1_{\bC_m}(X_t)-\mu_b(\bC_m)\big)\,
\1_{\{\tau_m\leq t-\Delta\}}\bigr]\,dt.
\end{align*}
Using the Markov property of the process $(X_t)$ and the spectral
gap inequality from Assumption \ref{A1}, we infer that
\begin{align*}
\Ex_b\bigl[\big(\1_{\bC_m}(X_t)&-\mu_b(\bC_m)\big)\,
\1_{\{\tau_m\leq
t-\Delta\}}\bigr]\\
&=\Ex_b\bigl[P_{b,\Delta}(\1_{\bC_m}-\mu_b(\bC_m))(X_{t-\Delta})\,
\1_{\{\tau_m\leq t-\Delta\}}\bigr]\\
&\leq
\sqrt{\Ex_b\bigl[\big(P_{b,\Delta}(\1_{\bC_m}-\mu_b(\bC_m))(X_{t-\Delta})\big)^2\bigr]}
\\
&=\|P_{b,\Delta}\1_{\bC_m}-\mu_b(\bC_m)\|_{\mu_b}\leq
e^{-\Delta\rho}\sqrt{\mu_b(\bC_m)}.
\end{align*}
This inequality completes the proof of the lemma. \qed
\end{proof}

\begin{lemma}
We have uniformly over $m=1,\ldots,M$:
$$
\Pb_b(\tau_m\le t)\lesssim \frac{t^2\eta^2+t\psi_d^2(h^d)}{(T-t)^2}.
$$
\end{lemma}

\begin{proof}
Note that $M_t:=\mathcal J_m^{-1/2}\hat{\eta}_{mj}(t)\in\RR^K$ is a
martingale with quadratic variation matrix $\langle
M\rangle_t=\mathcal J_m^{-1/2}\hat{\mathcal J}_m(t)\mathcal
J_m^{-1/2}$. We obtain that  $\Ex_b[\langle M\rangle_t]=t I_K$ with
the $K\times K$-unit matrix $I_K$ and
\begin{align*}
\Pb_b(\tau_m\le t)&=\Pb_b(\|\langle M\rangle_t\|\ge T)=
\Pb_b(\|\langle M\rangle_t-tI_K\|\ge T-t)\\
&\le \frac{\Ex_b[\|\langle M\rangle_t-tI_K\|^2]}{(T-t)^2}\ .
\end{align*}
Let $J_h\in\RR^{K\times K}$ be the diagonal matrix with
$J_{h,ii}=h^{|\alpha(i)|}$, $i=1,\ldots,K$, then
\begin{align*}
\|\langle M\rangle_t-tI_K\|&=\|\mathcal J_m^{-1/2}(\hat{\mathcal
J}_m(t)-t\mathcal
J_m)\mathcal J_m^{-1/2}\|\\
&\le\|\mathcal J_m^{-1/2}J_h\|^2\|J_h^{-1}(\hat{\mathcal
J}_m(t)-t\mathcal J_m)J_h^{-1}\|.
\end{align*}
Simple algebra shows that $\|\mathcal
J_m^{-1/2}J_h\|^2=\|(J_h^{-1}\mathcal J_m J_h^{-1})^{-1}\|$,
$J_h^{-1}=J_{h^{-1}}$ and
$$
J_h^{-1}\mathcal J_m
J_h^{-1}=h^d\int_{[0,1]^d}\mv(u)\mv(u)^T\mu_{b^{\circ}}(a_m+uh)\,du.
$$
This matrix is strictly positive definite and
$\|h^{-d}\mu_{b^{\circ}}(a_m)^{-1}J_h^{-1}\mathcal J_m
J_h^{-1}-\mV\|$ tends to zero as $h\to 0$. Hence, by the continuity
of the matrix inversion we obtain for $h$ small enough
$$
\|h^{d}\mu_{b^{\circ}}(a_m)J_h\mathcal J_m^{-1} J_h\|\leq
2\|\mV^{-1}\|.
$$
We conclude that $\|\mathcal J_m^{-1/2}J_h\|^2\lesssim
\mu_{b^{\circ}}(\bC_m)^{-1}$. Set now $H_t=J_h^{-1}(\hat{\mathcal
J}_m(t)-t\mathcal J_m)J_h^{-1}$. It is easily checked that
\begin{align*}
H_t=&\int_0^t\1_{\bC_m}(X_s)\mv\Big(\frac{X_s-a_m}{h}\Big)\mv\Big(\frac{X_s-a_m}{h}\Big)^T\!ds\\
&-
t\int_{\bC_m}\mv\Big(\frac{x-a_m}{h}\Big)\mv\Big(\frac{x-a_m}{h}\Big)^T\!\mu_{b^{\circ}}(x)\,dx.
\end{align*}
Each entry $H_{t,ij}$ can be written as $\int_0^t
f(X_s)\,ds-t\int_{\bC_m}f(x)\mu_{b^{\circ}}(x)\,dx$, where $f$ is a
function bounded by $1$ and supported by $\bC_m$. Thus, a
bias-variance decomposition combined with Proposition~\ref{adfunc}
yields
$$
\Ex_b[H_{t,ij}^2]\lesssim
t^2\bigg(\int_{\bC_m}|\mu_b(x)-\mu_{b^{\circ}}(x)|\,dx\bigg)^2+th^d\psi_d^2(h^d)\mu_b(\bC_m).
$$
Since in view of Remark~\ref{lowmu} $\mu_b(\bC_m)$ and
$\mu_{b^{\circ}}(\bC_m)$ are both of order $h^d$ and all norms in
$\RR^{K\times K}$ are equivalent, we arrive at the desired estimate.
\qed\end{proof}

Using the last lemma we obtain
\begin{align*}
\int_0^{T-\Delta}\!\!\Pb_b(\tau_m\leq
t)\,dt&\lesssim \int_0^T \min\Big(1,\frac{t^2\eta^2}{(T-t)^2}+\frac{\psi_d^2(h^d)t}{(T-t)^{2}}\Big)\,dt\\
&\leq \int_0^T\!\! \min\Big(1,\frac{t^2\eta^2}{(T-t)^2}\Big)\,dt+
\int_0^T\!\! \min\Big(1,\frac{\psi_d^2(h^d)t}{(T-t)^{2}}\Big)\,dt.
\end{align*}

Setting $c_T=T^{-1/2}\psi_d(h^d)$, we get
\begin{align*}
\int_0^T \min\Big(1,\frac{\psi_d^2(h^d)t}{(T-t)^{2}}\Big)\,dt&=T\int_0^1\min(1,c_T^2(1-v)v^{-2})\,dv\\
&\le
T\int_0^{c_T}1\,dv+T\int_{c_T}^\infty
c_T^2v^{-2}\,dv\\
&=2Tc_T=2T^{1/2}\psi_d(h^d).
\end{align*}
In the same way we obtain
$\int_0^T\min\big(1,{t^2\eta^2}/{(T-t)^2}\big)\,dt\leq 2T\eta$.
Substituting all estimates into \eqref{intdecomp} and
\eqref{T3decomp}, we obtain
$$
\mathcal T_3(\vartheta)\lesssim \|\bar b-\bar
b^{\circ}\|_\infty^2\big(\Delta+T\eta+\psi_d(h^d)\sqrt{T}+Th^{-d/2}e^{-\Delta\rho}\big).
$$
Thus choosing $\Delta(T)=\psi_d(h^d)\sqrt{T}$ we get
\[\mathcal
T_3(\vartheta)\lesssim \|\bar b-\bar
b^{\circ}\|_\infty^2\,(T\eta+\psi_d(h^d)\sqrt{T}),
\] \
provided that $h=h(T)$ tends to zero as $T\to\infty$.

\end{document}